\documentclass[12pt,a4paper,reqno]{amsart}%
\usepackage{amsmath}
\usepackage{amsfonts}
\usepackage{amssymb}
\usepackage{graphicx}
\usepackage[colorlinks]{hyperref}
\usepackage[a4paper,left=3cm,right=3cm]{geometry}%
\providecommand{\U}[1]{\protect\rule{.1in}{.1in}}
\allowdisplaybreaks[4]
\newtheorem{theorem}{Theorem}

\newtheorem{definition}[theorem]{Definition}

\newtheorem{lemma}[theorem]{Lemma}

\newtheorem{proposition}[theorem]{Proposition}

\begin{document}
\title{Irregularities of distribution on two point homogeneous spaces}
\author[L. Brandolini]{Luca Brandolini}
\address[L. Brandolini]{Dipartimento di Ingegneria Gestionale, dell'Informazione e
della Produzione, Universit\`a degli Studi di Bergamo, Viale Marconi 5,
Dalmine BG, Italy}
\email{luca.brandolini@unibg.it}
\author[B. Gariboldi]{Bianca Gariboldi}
\address[B. Gariboldi]{Dipartimento di Matematica, Universit\`a degli Studi di Milano
Bicocca, via Roberto Cozzi 55, Milano, Italy}
\email{biancamaria.gariboldi@guest.unibg.it}
\author[G. Gigante]{Giacomo Gigante}
\address[G. Gigante]{Dipartimento di Ingegneria Gestionale, dell'Informazione e della
Produzione, Universit\`a degli Studi di Bergamo, Viale Marconi 5, Dalmine BG, Italy}
\email{giacomo.gigante@unibg.it}

\begin{abstract}
We study the irregularities of distribution on two-point homogeneous spaces.
Our main result is the following: let $d$ be the real dimension of a two point
homogeneous space $\mathcal{M}$, let $\left(  \{ a_{j}\} _{j=1}^{N},\{ x_{j}\}
_{j=1}^{N}\right)  $ be a system of positive weights and points on
$\mathcal{M}$ and let
\[
D_{r}( x) =\sum_{j=1}^{N}a_{j}\chi_{B_{r}(x)}(x_{j})-\mu(B_{r}(x))
\]
be the discrepancy associated with the ball $B_{r}( x) $. Then, if
$d\not \equiv 1(\operatorname{mod}4)$, for any radius $0<r<\pi/2$, we obtain
the sharp estimate%
\[
\int_{\mathcal{M}}\left(  \left\vert D_{r}( x) \right\vert ^{2}+\left\vert
D_{2r}( x) \right\vert ^{2}\right)  d\mu( x) \geqslant cN^{-1-\frac{1}{d}}.
\]

\end{abstract}
\subjclass[2010]{11K38, 43A85}
\keywords{Irregularities of distribution, Two-point homogeneous spaces, Discrepancy}
\maketitle

\section{Introduction}

A $d$-dimensional Riemannian manifold $\mathcal{M}$ with distance $\rho$ is
said to be a two-point homogeneus space if given four points $x_{1}%
,x_{2},y_{1},y_{2}\in\mathcal{M}$ such that $\rho(x_{1},y_{1})=\rho
(x_{2},y_{2})$, then there exists an isometry $g$ of $\mathcal{M}$ such that
$gx_{1}=x_{2}$ and $gy_{1}=y_{2}$. Wang in \cite{Wang} has completely
characterized compact connected two-point homogeneous spaces. More precisely
if we assume that $\rho$ is normalized so that $\mathrm{diam}(\mathcal{M}%
)=\pi$, then $\mathcal{M}$ is isometric to one of the following compact rank
$1$ symmetric spaces:

\begin{itemize}
\item[(i)] the Euclidean sphere $S^{d}=SO(d+1) /SO(d)\times\{1\}$,
$d\geqslant1$;

\item[(ii)] the real projective space $P^{n}( \mathbb{R}) =O(n+1) /O( n)
\times O( 1) $, $n\geqslant2$;

\item[(iii)] the complex projective space $P^{n}( \mathbb{C}) =U(n+1) /U( n)
\times U( 1) $, $n\geqslant2$;

\item[(iv)] the quaternionic projective space $P^{n}( \mathbb{H})=Sp(n+1)/Sp(
n) \times Sp( 1) $, $n\geqslant2$;

\item[(v)] the octonionic projective plane $P^{2}( \mathbb{O}) $.
\end{itemize}

In the following we will assume that $\mathcal{M}$ is one of the above
symmetric spaces and that $d$ is its real dimension. In particular the real
dimension of $P^{n}(\mathbb{F})$ is $d=nd_{0}$, where $d_{0}=1,2,4,8$
according to the real dimension of $\mathbb{F}=\mathbb{R}$, $\mathbb{C}$,
$\mathbb{H}$ and $\mathbb{O}$ respectively. In the case of $S^{d}$ it will be
convenient to set $d_{0}=d$. See \cite[pp. 176-178]{Gangolli}, see also
\cite{H}, \cite{Skriganov} and \cite{Wolf}. In the following, to keep notation
simple, we will use
\[
a =\frac{d-2}{2}, \ \quad\ b =\frac{d_{0}-2}{2}.
\]

Let $\mu$ be the Riemannian measure on $\mathcal{M}$ normalized so that
$\mu(\mathcal{M})=1$ and let $B_{r}(x)=\{y\in\mathcal{M}:\rho(x,y)<r\}$.

For a given set of points $\left\{  x_{j}\right\}  _{j=1}^{N}\subset
\mathcal{M}$ and positive weights $\left\{  a_{j}\right\}  _{j=1}^{N}$ such
that $a_{1}+a_{2}+\cdots+a_{N}=1$ we define the discrepancy of the ball
$B_{r}(x)$ by
\begin{equation}
D_{r}(x)=\sum_{j=1}^{N}a_{j}\chi_{B_{r}(x)}(x_{j})-\mu(B_{r}(x))
\label{Discrepanza M}%
\end{equation}

This quantity compares the Riemannian measure $\mu$ with the discrete measure
$\sum_{j=1}^{N}a_{j}\delta_{x_{j}}$ testing these two measures on the ball
$B_{r}(x)$. It is known that, no matter how well distributed the points
$\{x_{j}\}_{j=1}^{N}$ are, there are balls for which the discrepancy
$D_{r}(x)$ cannot be too small. See \cite{BC} for an introduction to this
subject. For example M. Skriganov in \cite[Theorem 2.2]{Skriganov} has proved,
in the case of equal weights, that if $\eta$ is a positive, locally integrable
function on $(0,\pi)$ such that%
\[
\int_{0}^{\pi}\eta(r)\left(  \sin\frac{1}{2}r\right)  ^{d-1}\left(  \cos
\frac{1}{2}r\right)  ^{d_{0}-1}dr<+\infty,
\]
then there exists $c>0$ such that for every distribution of $N$ points%
\[
\int_{0}^{\pi}\int_{\mathcal{M}}\left\vert D_{r}(x)\right\vert ^{2}d\mu
(x)\eta(r)dr\geqslant cN^{-1-\frac{1}{d}}.
\]
A very general non-constructive result guarantees the existence of point
distributions that exhibit the optimal decay $N^{-1-1/d}$ (see e.g.
\cite[Corollary 8.2]{BCCGT}). In fact, Skriganov \cite[Corollary
2.1]{Skriganov} has proved that for well distributed optimal cubature formulas
one has
\begin{equation}
\int_{0}^{\pi}\int_{\mathcal{M}}\left\vert D_{r}(x)\right\vert ^{2}d\mu
(x)\eta(r)dr\leqslant cN^{-1-\frac{1}{d}} \label{stima cubatura alto}%
\end{equation}
(see \cite{GG}, or the remarks that follow the statement of Theorem
\ref{Thm cubature alto} here, for a proof of the existence of such cubature formulas).

Our first result is the following.

\begin{theorem}
\label{Thm 1}Let $d\not \equiv 1 \ ( \operatorname{mod}4) $. Then, for all
$0<r<\pi/2$, there is a constant $C>0$ such that for every set of points
$\left\{  x_{j}\right\}  _{j=1}^{N}\subset\mathcal{M}$ and positive weights
$\left\{  a_{j}\right\}  _{j=1}^{N}$ such that $a_{1}+a_{2}+\cdots+ a_{N}=1$
we have
\[
\int_{\mathcal{M}}\left\vert D_{r}( x) \right\vert ^{2}d\mu( x) +\int
_{\mathcal{M}}\left\vert D_{2r}( x) \right\vert ^{2}d\mu( x) \geqslant
CN^{-1-\frac{1}{d}}.
\]

\end{theorem}

Again this result is optimal in view of Corollary 8.2 in \cite{BCCGT}. Also,
the same technique used by Skriganov to prove the estimate from above for
cubature formulas (\ref{stima cubatura alto}) actually gives a stronger,
uniform estimate in the radius $r$,%
\begin{equation}
\int_{\mathcal{M}}\left\vert D_{r}(x)\right\vert ^{2}d\mu(x)\leqslant
cN^{-1-\frac{1}{d}}. \label{dallalto}%
\end{equation}
See Theorem \ref{Thm cubature alto} for a precise statement.

When $d\equiv1\,\left(  \operatorname{mod}4\right)  $ our technique fails. In
similar settings it is known that the discrepancy can be actually a bit
smaller than the expected value $N^{-1-\frac{1}{d}}$. See \cite[Theorem
3.1]{PS}. See also \cite{BCGT}, \cite{KSS} and \cite{PSido}.

Our final result is an estimate for the discrepancy of balls of a single fixed
radius that holds in every dimension.

\begin{theorem}
\label{Thm 2}For every $\varepsilon>0$ and for almost every $0<r<\pi$ there is
a constant $C>0$ such that for every set of points $\left\{  x_{j}\right\}
_{j=1}^{N}\subset\mathcal{M}$ and positive weights $\left\{  a_{j}\right\}
_{j=1}^{N}$ such that $a_{1}+a_{2}+\cdots+ a_{N}=1$ we have
\[
\int_{\mathcal{M}}\left\vert D_{r}(x)\right\vert ^{2}d\mu(x)\geqslant
CN^{-1-\frac{3}{d}-\varepsilon}.
\]

\end{theorem}

The paper is organized as follows. In section 2 we collect some basic facts
about two-point homogeneous spaces. In section 3 we prove Theorem \ref{Thm 1}
and \ref{Thm 2}. In section 4 we show that suitable cubature formulas give the
optimal discrepancy (\ref{dallalto}).

\section{Two-point homogeneous spaces}

If $o$ is a fixed point in $\mathcal{M}$, then $\mathcal{M}$ can be identified
with the homogeneous space $G/K$, where $G$ is the group of isometries of
$\mathcal{M}$ and $K$ is the stabilizer of $o$. We will also identify
functions $F(x)$ on $\mathcal{M}$ with right $K$-invariant functions $f(g)$ on
$G$ by setting $f(g)=F(x)$ when $go=x$. If
$\mu$ is the Riemannian measure on $\mathcal{M}$ normalized so that
$\mu(\mathcal{M})=1$,
then $\mu$ is invariant under the action of $G$,
in other words, for every $g\in G$,
\[
\int_{\mathcal{M}}F(gx)d\mu(x)=\int_{\mathcal{M}}F(x)d\mu(x)
\]

\begin{definition}
A function $F$ on $\mathcal{M}$ is a zonal function if for every
$x\in\mathcal{M}$ and every $k\in K$ we have $F( kx) =F( x) $.
\end{definition}

\begin{lemma}
\label{Lemma misura}Let $F$ be a zonal function. Then $F( x) $ depends only on
$\rho( x,o) $. Furthermore, defining $F_{0}$ so that $F( x) =F_{0}( \rho( x,o)
) $ we have%
\begin{equation}
\int_{\mathcal{M}}F( x) d\mu( x) =\int_{0}^{\pi}F_{0}( r) A( r) dr,
\label{IntegraleZonale}%
\end{equation}
where%
\[
A( r) =c( a,b) \left(  \sin\frac{1}{2}r\right)  ^{2a+1}\left(  \cos\frac{1}%
{2}r\right)  ^{2b+1}%
\]
and
\[
c( a,b) =\left(  \int_{0}^{\pi}\left(  \sin\frac{1}{2}r\right)  ^{2a+1}\left(
\cos\frac{1}{2}r\right)  ^{2b+1}dr\right)  ^{-1} =\frac{\Gamma(a+b+2)}%
{\Gamma(a+1)\Gamma(b+1)}.
\]
In particular if%
\[
B_{r}( x) =\{ y\in\mathcal{M}:\rho( x,y) <r\} ,
\]
there exist two constants $c_{1}( d,d_{0}) $ and $c_{2}( d,d_{0}) $ such that
for every $r\in[ 0,\pi] $%
\begin{equation}
c_{1}( d,d_{0}) r^{d}\leqslant\mu( B_{r}( x) ) \leqslant c_{2}( d,d_{0})
r^{d}. \label{Volume palle}%
\end{equation}
Moreover if $k\geqslant1$, then a ball of radius $kr$ can be covered by
$\frac{c_{2}( d,d_{0}) }{c_{1}( d,d_{0}) }( 2(k+1)) ^{d}$ balls of radius $r$.
\end{lemma}

\begin{proof}
Let $x,y\in\mathcal{M}$ such that $\rho( x,o) =\rho( y,o) $. Since
$\mathcal{M}$ is two-point homogeneous there exists $g\in G$ such that $gx=y$
and $go=o$. Thus, $g\in K$ and $F( y) =F( gx) =F( x) $. Equation
(\ref{IntegraleZonale}) follows from (4.17) in \cite{Gangolli}. The estimate
(\ref{Volume palle}) is an immediate consequence of (\ref{IntegraleZonale}).
The last assertion is an estimate of the maximum number of disjoint balls of
radius $r/2$ that can fit in a ball of radius $(k+1)r$.
\end{proof}

Let $\Delta$ be the Laplace-Beltrami operator on $\mathcal{M}$, let
$\lambda_{0},\lambda_{1},\ldots,$ be the distinct eigenvalues of $-\Delta$
arranged in increasing order, let $\mathcal{H}_{m}$ be the eigenspace
associated with the eigenvalue $\lambda_{m}$, and let $d_{m}$ its dimension.
It is well known that%
\begin{equation}
L^{2}(\mathcal{M})=%
{\displaystyle\bigoplus_{m=0}^{+\infty}}
\mathcal{H}_{m}. \label{Decomp L2}%
\end{equation}
If $F(x)=F_{0}(\rho(x,o))$ is a zonal function on $\mathcal{M}$, then%
\begin{equation}
\Delta F(x)=\left.  \frac{1}{A(t)}\frac{d}{dt}\left(  A(t)\frac{d}{dt}%
F_{0}(t)\right)  \right\vert _{t=\rho(x,o)} \label{Radial Laplace Beltrami}%
\end{equation}
(see (4.16) in \cite{Gangolli}).

\begin{definition}
The zonal spherical function of degree $m\in\mathbb{N}$ with pole
$x\in\mathcal{M}$ is the unique function $Z_{x}^{m}\in\mathcal{H}_{m}$, given
by the Riesz representation theorem, such that for every $Y\in\mathcal{H}_{m}$%
\[
Y( x) =\int_{\mathcal{M}}Y( y) Z_{x}^{m}( y) d\mu( y) .
\]

\end{definition}

The next lemma summarizes the main properties of zonal functions and it is
essentially taken from \cite{SW} where the case $\mathcal{M}=S^{d}$ is
discussed in detail.

\begin{lemma}
\begin{itemize}

\item[i)] If $Y_{m}^{1},\ldots,Y_{m}^{d_{m}}$ is an orthonormal basis of
$\mathcal{H}_{m}\subset L^{2}( \mathcal{M}) $, then%
\[
Z_{x}^{m}( y) =\sum_{\ell=1}^{d_{m}}\overline{Y_{m}^{\ell}( x) }Y_{m}^{\ell}(
y) .
\]

\item[ii)] $Z_{x}^{m}$ is real valued and $Z_{x}^{m}( y) =Z_{y}^{m}( x) .$

\item[iii)] If $g\in G$, then $Z_{gx}^{m}( gy) =Z_{x}^{m}( y) .$

\item[iv)] $Z_{x}^{m}( x) =\Vert Z_{x}^{m}\Vert_{2}^{2}=d_{m}.$

\item[v)] $Z_{o}^{m}( x) $ is a zonal function and
\begin{equation}
Z_{o}^{m}( x) =\frac{d_{m}}{P_{m}^{a,b}( 1) }P_{m}^{a,b}( \cos( \rho( x,o) ) )
\label{Jacobi}%
\end{equation}
where $P_{m}^{a,b}$ are the Jacobi polynomials.

\item[vi)] $\left\{  d_{m}^{-1/2}Z_{o}^{m}\right\}  _{m=0}^{+\infty}$ is an
orthonormal basis of the subspace of $L^{2}(\mathcal{M})$ of zonal functions.

\item[vii)] $\lambda_{m}=m(m+a+b+1)$.

\item[viii)] $d_{m}\approx m^{d-1}$.
\end{itemize}
\end{lemma}

\begin{proof}
i) By the defining property of $Z_{x}^{m}$%
\[
Z_{x}^{m}(y)=\sum_{\ell=1}^{d_{m}}\int_{\mathcal{M}}Z_{x}^{m}(z)\overline
{Y_{m}^{\ell}(z)}d\mu(z)Y_{m}^{\ell}(y)=\sum_{\ell=1}^{d_{m}}\overline
{Y_{m}^{\ell}(x)}Y_{m}^{\ell}(y).
\]
ii) Since the basis can be taken to be of real valued functions, then by point
i) $Z_{x}^{m}$ is real valued.\newline iii) Let $Y\in\mathcal{H}%
_{m}(\mathcal{M})$. Since for every $g$, $Y(g^{-1}y)\in\mathcal{H}%
_{m}(\mathcal{M})$, then
\[
\int_{\mathcal{M}}Y(y)Z_{gx}^{m}(gy)d\mu(y)=\int_{\mathcal{M}}Y(g^{-1}%
y)Z_{gx}^{m}(y)d\mu(y)=Y(g^{-1}(gx))=Y(x),
\]
and the uniqueness of $Z_{x}^{m}$ shows that $Z_{gx}^{m}(gy)=Z_{x}^{m}%
(y)$.\newline iv) Observe that%
\[
\Vert Z_{x}^{m}\Vert_{2}^{2}=\int_{\mathcal{M}}Z_{x}^{m}(y)\overline{Z_{x}%
^{m}(y)}d\mu(y)=\overline{Z_{x}^{m}(x)}=\sum_{\ell=1}^{d_{m}}|Y_{m}^{\ell
}(x)|^{2}.
\]
Since by iii) this quantity is constant in $x$, it equals its integral over
$\mathcal{M}$, which is exactly $d_{m}$.\newline v) By point iii) we have%
\[
Z_{o}(kx)=Z_{k^{-1}o}(x)=Z_{o}(x).
\]
Equation (\ref{Jacobi}) follows from the fact that the radial part of the
Laplace-Beltrami operator (\ref{Radial Laplace Beltrami}) coincides with the
Jacobi operator.\newline vi) follows from the fact that the Jacobi
polynomials
\[
\left\{  P_{m}^{a,b}(\cos(t))\right\}  _{m=0}^{+\infty}%
\]
form an orthogonal basis of $L^{2}((0,\pi),A(t)dt)$. See Chapter 4 in
\cite{szego}. \newline vii) follows from \cite[Theorem 4.2.2]{szego}. \newline
viii) follows from iv), v) and vi) computing explicitly the $L^{2}$ norm of
$d_{n}^{-1/2}Z_{o}^{m}$. For the computation one needs (4.1.1) in \cite{szego}
to evaluate $P_{m}^{a,b}(1)$ and (4.3.3) in \cite{szego} to evaluate the
$L^{2}$ norm of $P_{m}^{a,b}(x)$.
\end{proof}

\section{Proof of the main results}

For the proof of our main results we need some lemmas. In the next one we
obtain an expression for the discrepancy that separates the contribution of
the point distribution (and weights) from the contribution of the geometry of
the set $B_{r}( x) $.

\begin{lemma}
\label{Lemma Quadrato Discrep}Let $D_{r}( x) $ be as in (\ref{Discrepanza M}).
Then%
\[
\int_{\mathcal{M}}\left\vert D_{r}( x) \right\vert ^{2}d\mu( x) =\sum
_{m=1}^{+\infty}\sum_{\ell=1}^{d_{m}}\left\vert \sum_{j=1}^{N}a_{j}Y_{m}%
^{\ell}( x_{j}) \right\vert ^{2}d_{m}^{-2}\left\vert \int_{B_{r}( o) }%
Z_{o}^{m}( y) d\mu( y) \right\vert ^{2}.
\]

\end{lemma}

\begin{proof}
By (\ref{Decomp L2}) we have%
\[
\int_{\mathcal{M}}\left\vert D_{r}(x)\right\vert ^{2}d\mu(x)=\sum
_{m=0}^{+\infty}\int_{\mathcal{M}}|\mathbb{P}_{m}D_{r}(x)|^{2}d\mu(x)
\]
(here $\mathbb{P}_{m}$ denotes the orthogonal projection of $L^{2}(
\mathcal{M}) $ onto $\mathcal{H}_{m}$). Let $g\in G$ such that $gx_{j}=o$.
Since $\chi_{B_{r}(o)}$ is a zonal function, then%
\begin{align*}
\chi_{B_{r}(x_{j})}(x)  &  =\chi_{B_{r}(o)}(gx)=\sum_{m=0}^{+\infty}d_{m}%
^{-1}\int_{B_{r}(o)}Z_{o}^{m}(y)d\mu(y)Z_{o}^{m}(gx)\\
&  =\sum_{m=0}^{+\infty}d_{m}^{-1}\int_{B_{r}(o)}Z_{o}^{m}(y)d\mu(y)Z_{x_{j}%
}^{m}(x).
\end{align*}
Since $Z_{x_{j}}^{m}(x)\in\mathcal{H}_{m}$, we have%
\[
\mathbb{P}_{m}\chi_{B_{r}(x_{j})}(x)=d_{m}^{-1}\int_{B_{r}(o)}Z_{o}^{m}%
(y)d\mu(y)Z_{x_{j}}^{m}(x).
\]
Overall we obtain%
\[
\mathbb{P}_{m}D_{r}(x)=\sum_{j=1}^{N}a_{j}d_{m}^{-1}\int_{B_{r}(o)}Z_{o}%
^{m}(y)d\mu(y)Z_{x_{j}}^{m}(x)-\delta_{0}(m)\mu(B_{r}(o)).
\]
In particular $\mathbb{P}_{0}D_{r}(x)=0$ and for $m>0$%
\begin{align*}
\mathbb{P}_{m}D_{r}(x)  &  =d_{m}^{-1}\sum_{j=1}^{N}a_{j}\int_{B_{r}(o)}%
Z_{o}^{m}(y)d\mu(y)Z_{x_{j}}^{m}(x)\\
&  =d_{m}^{-1}\sum_{\ell=1}^{d_{m}}\left(  \sum_{j=1}^{N}a_{j}\overline
{Y_{m}^{\ell}(x_{j})}\right)  \int_{B_{r}(o)}Z_{o}^{m}(y)d\mu(y)Y_{m}^{\ell
}(x)
\end{align*}
Finally%
\[
\int_{\mathcal{M}}|D_{r}(x)|^{2}d\mu(x)=\sum_{m=1}^{+\infty}\sum_{\ell
=1}^{d_{m}}\left\vert \sum_{j=1}^{N}a_{j}Y_{m}^{\ell}(x_{j})\right\vert
^{2}d_{m}^{-2}\left\vert \int_{B_{r}(o)}Z_{o}^{m}(y)d\mu(y)\right\vert ^{2}.
\]

\end{proof}

In the next results we estimate the quantities%
\[
\sum_{\ell=1}^{d_{m}}\left\vert \sum_{j=1}^{N}a_{j}Y_{m}^{\ell}(x_{j}%
)\right\vert ^{2} \quad\text{and } \quad d_{m}^{-2}\left\vert \int_{B_{r}%
(o)}Z_{o}^{m}(y)d\mu(y)\right\vert ^{2}.
\]

\begin{proposition}
\label{Prop Cassels}Let $m_{0}>0$, then there exist $C_{0},C_{1}>0$ such that
for every $X\geqslant m_{0}$ we have%
\[
\sum_{m=m_{0}}^{X}\sum_{\ell=1}^{d_{m}}\left\vert \sum_{j=1}^{N}a_{j}%
Y_{m}^{\ell}(x_{j})\right\vert ^{2}\geqslant C_{1}\sum_{j=1}^{N}a_{j}^{2}%
X^{d}-C_{0}.
\]

\end{proposition}

\begin{proof}
By the Cassels-Montgomery inequality for manifolds (see \cite[Theorem 1]{BGG})
along with the fact that $\sum_{m=0}^{X}d_{m}\approx\sum_{m=0}^{X}%
m^{d-1}\approx X^{d}$, we have%
\[
\sum_{m=0}^{X}\sum_{\ell=1}^{d_{m}}\left\vert \sum_{j=1}^{N}a_{j}Y_{m}^{\ell
}(x_{j})\right\vert ^{2}\geqslant C_{1}\sum_{j=1}^{N}a_{j}^{2}X^{d}.
\]
Since%
\begin{align*}
\sum_{m=0}^{m_{0}-1}\sum_{\ell=1}^{d_{m}}\left\vert \sum_{j=1}^{N}a_{j}%
Y_{m}^{\ell}(x_{j})\right\vert ^{2}  &  \leqslant\max_{\substack{0\leqslant
m\leqslant m_{0}-1\\\ell=1,\ldots,d_{m}}}\Vert Y_{m}^{\ell}\Vert_{\infty}%
^{2}\sum_{m=0}^{m_{0}-1}\sum_{\ell=1}^{d_{m}}\left\vert \sum_{j=1}^{N}%
a_{j}\right\vert ^{2}\\
&  =\max_{\substack{0\leqslant m\leqslant m_{0}-1\\\ell=1,\ldots,d_{m}}}\Vert
Y_{m}^{\ell}\Vert_{\infty}^{2}\sum_{m=0}^{m_{0}-1}d_{m}=C_{0},
\end{align*}
we have%
\[
\sum_{m=m_{0}}^{X}\sum_{\ell=1}^{d_{m}}\left\vert \sum_{j=1}^{N}a_{j}%
Y_{m}^{\ell}(x_{j})\right\vert ^{2}\geqslant C_{1}\sum_{j=1}^{N}a_{j}^{2}%
X^{d}-C_{0}.
\]

\end{proof}

\begin{lemma}
\label{Trasf bolla}For any $0\leqslant r\leqslant\pi$ and for any
$m\geqslant1$ we have%
\[
\int_{B_{r}(o)}Z_{o}^{m}(x)d\mu(x)=\frac{c( a,b) d_{m}}{mP_{m}^{a,b}%
(1)}P_{m-1}^{a+1,b+1}(\cos r)\left(  \sin\left(  \frac{r}{2}\right)  \right)
^{2a+2}\left(  \cos\left(  \frac{r}{2}\right)  \right)  ^{2b+2}.
\]

\end{lemma}

\begin{proof}
By (\ref{IntegraleZonale}), (\ref{Jacobi}), and a change of variable,%
\begin{align*}
\int_{B_{r}(o)}Z_{o}^{m}(x)d\mu(x)  &  =\frac{d_{m}}{P_{m}^{a,b}(1)}\int
_{0}^{r}P_{m}^{a,b}(\cos t)A(t)dt\\
&  =\frac{c(a,b)d_{m}}{P_{m}^{a,b}(1)}\int_{0}^{r}P_{m}^{a,b}(\cos t)\left(
\sin\frac{t}{2}\right)  ^{2a+1}\left(  \cos\frac{t}{2}\right)  ^{2b+1}dt\\
&  =\frac{c(a,b)d_{m}}{2P_{m}^{a,b}(1)}\int_{0}^{r}P_{m}^{a,b}(\cos t)\left(
\frac{1-\cos t}{2}\right)  ^{a}\left(  \frac{1+\cos t}{2}\right)  ^{b}\sin
tdt\\
&  =\frac{c(a,b)d_{m}}{2^{a+b+1}P_{m}^{a,b}(1)}\int_{\cos r}^{1}P_{m}%
^{a,b}(x)(1-x)^{a}(1+x)^{b}dx.
\end{align*}
By Rodrigues' formula (see \cite[(4.3.1)]{szego})%
\begin{equation}
\label{rodrigues}P_{m}^{a,b}(x)(1-x)^{a}(1+x)^{b}=-\frac{1}{2m}\frac{d}%
{dx}\left(  P_{m-1}^{a+1,b+1}(x)(1-x)^{a+1}(1+x)^{b+1}\right)  ,
\end{equation}
thus%
\begin{align*}
\int_{B_{r}(o)}Z_{o}^{m}(x)d\mu(x)  &  =\frac{c(a,b)d_{m}}{2m2^{a+b+1}%
P_{m}^{a,b}(1)}P_{m-1}^{a+1,b+1}(\cos r)(1-\cos r)^{a+1}(1+\cos r)^{b+1}\\
&  =\frac{c(a,b)d_{m}}{mP_{m}^{a,b}(1)}P_{m-1}^{a+1,b+1}(\cos r)\left(
\sin\left(  \frac{r}{2}\right)  \right)  ^{2a+2}\left(  \cos\left(  \frac
{r}{2}\right)  \right)  ^{2b+2}.
\end{align*}

\end{proof}

In the following, to keep notation simple we set
\[
M=m+\frac{a+b+1}{2}.
\]

\begin{lemma}
\label{coeff_bolla}For all $\varepsilon>0$ and for all integers $m\geqslant1$,%
\begin{align*}
&  \int_{B_{r}(o)}Z_{o}^{m}(y)d\mu(y)\\
=  &  \ c(a,b)d_{m}\Gamma(a+1)\left(  \sin\left(  \frac{r}{2}\right)  \right)
^{a+1}\left(  \cos\left(  \frac{r}{2}\right)  \right)  ^{b+1}M^{-a-1}\left(
\frac{r}{\sin r}\right)  ^{\frac{1}{2}}J_{a+1}( Mr)\\
&  +d_{m}O( m^{-\frac{5}{2}-a}),
\end{align*}
uniformly in $r\in[0,\pi-\varepsilon]$. Here $J_{a+1}$ is the Bessel function
of first type and order $a+1$.
\end{lemma}

\begin{proof}
The following asymptotic expansion of the Jacobi polynomials is well known,
(see \cite[Theorem 8.21.12]{szego}),%
\begin{align}
&  \left(  \sin\frac{r}{2}\right)  ^{a+1}\left(  \cos\frac{r}{2}\right)
^{b+1}P_{m-1}^{a+1,b+1}(\cos r)\label{Jacobi Bessel}\\
=  &  M^{-a-1}\frac{\Gamma(m+a+1)}{(m-1)!}\left(  \frac{r}{\sin r}\right)
^{\frac{1}{2}}J_{a+1}\left(  Mr\right)  +r^{\frac{1}{2}}O( m^{-\frac{3}{2}%
})\nonumber
\end{align}
Thus, by the previous lemma, recalling that (see (4.1.1) in \cite{szego})
\[
P_{m}^{a,b}(1)=\frac{\Gamma(m+a+1)}{\Gamma(m+1)\Gamma(a+1)},
\]
for every $r\in[0,\pi-\varepsilon]$ and $m\geqslant1$,
\begin{align*}
&  \int_{B_{r}(o)}Z_{o}^{m}(y)d\mu(y)\\
=  &  \ \frac{c(a,b)d_{m}}{mP_{m}^{a,b}(1)}P_{m-1}^{a+1,b+1}(\cos r)\left(
\sin\left(  \frac{r}{2}\right)  \right)  ^{2a+2}\left(  \cos\left(  \frac
{r}{2}\right)  \right)  ^{2b+2}\\
=  &  \ \frac{c(a,b)d_{m}}{m}\frac{\Gamma(m+1)\Gamma(a+1)}{\Gamma
(m+a+1)}\left(  \sin\left(  \frac{r}{2}\right)  \right)  ^{a+1}\left(
\cos\left(  \frac{r}{2}\right)  \right)  ^{b+1}\\
&  \times M^{-a-1}\frac{\Gamma(m+a+1)}{(m-1)!}\left(  \frac{r}{\sin r}\right)
^{1/2}J_{a+1}( Mr)\\
&  +\frac{c(a,b)}{m}\frac{d_{m}}{P_{m}^{a,b}(1)}\left(  \sin\left(  \frac
{r}{2}\right)  \right)  ^{a+1}\left(  \cos\left(  \frac{r}{2}\right)  \right)
^{b+1}r^{\frac{1}{2}}O(m^{-\frac{3}{2}})\\
=  &  \ c(a,b)d_{m}\Gamma(a+1)\left(  \sin\left(  \frac{r}{2}\right)  \right)
^{a+1}\left(  \cos\left(  \frac{r}{2}\right)  \right)  ^{b+1}M^{-a-1}\left(
\frac{r}{\sin r}\right)  ^{\frac{1}{2}} J_{a+1}\left(  Mr\right) \\
&  +d_{m}O(m^{-\frac{5}{2}-a}).
\end{align*}

\end{proof}

The expression obtained in Lemma \ref{Lemma Quadrato Discrep} shows that to
estimate the discrepancy from below requires an estimate from below for the
quantity
\[
\left\vert \int_{B_{r}(o)}Z_{o}^{m}(y)d\mu(y)\right\vert ^{2}.
\]
Lemma \ref{coeff_bolla} shows that this quantity may vanish due to the zeroes
of the Bessel functions. Our next lemma shows that at least when
$d\not \equiv 1\,( \operatorname{mod}4) $ one can overcome this obstruction
using two balls of different radii.

\begin{lemma}
\label{Lemma Due cerchi}Assume that $d\not \equiv 1\,( \operatorname{mod}4) $.
Then for all $0<r<\pi/2$, there exist positive constants $C$ and $m_{0}$ such
that for $m\geqslant m_{0}$,%
\[
\left\vert \int_{B_{r}(o)}Z_{o}^{m}(x)d\mu(x)\right\vert ^{2}+\left\vert
\int_{B_{2r}(o)}Z_{o}^{m}(x)d\mu(x)\right\vert ^{2}\geqslant Cd_{m}%
^{2}m^{-2a-3}.
\]

\end{lemma}

\begin{proof}
By the previous lemma, for all $m\geqslant1$%
\begin{align*}
&  \left\vert \int_{B_{r}(o)}Z_{o}^{m}(x)d\mu(x)\right\vert ^{2}+\left\vert
\int_{B_{2r}(o)}Z_{o}^{m}(x)d\mu(x)\right\vert ^{2}\\
\geqslant &  \ \left\vert c(a,b)d_{m}\Gamma(a+1)M^{-a-1}\right\vert
^{2}\left\{  \left\vert \left(  \sin\left(  \frac{r}{2}\right)  \right)
^{a+1}\left(  \cos\left(  \frac{r}{2}\right)  \right)  ^{b+1}\left(  \frac
{r}{\sin r}\right)  ^{\frac{1}{2}}J_{a+1}( Mr) \right\vert ^{2}\right. \\
&  + \left.  \left\vert (\sin r)^{a+1}(\cos r)^{b+1}\left(  \frac{2r}{\sin
2r}\right)  ^{\frac{1}{2}}J_{a+1}\left(  M2r\right)  \right\vert ^{2}\right\}
-d_{m}^{2}O( m^{-4-2a})\\
\geqslant &  \ Cd_{m}^{2}m^{-2a-2}\left\{  \vert J_{a+1}( Mr) \vert^{2}+\vert
J_{a+1}( 2Mr) \vert^{2}\right\}  -d_{m}^{2}O( m^{-4-2a}) .
\end{align*}
By the asymptotic expansion of the Bessel functions (see e.g. \cite[Lemma
3.11, Chapter 4]{SW}),%
\begin{equation}
J_{a+1}(w)=\sqrt{\frac{2}{\pi w}}\cos\left(  w-( a+1) \frac{\pi}{2}-\frac{\pi
}{4}\right)  +O( w^{-\frac{3}{2}}) , \label{asymp Bessel}%
\end{equation}
we have%
\begin{align*}
&  \left\vert J_{a+1}( Mr) \right\vert ^{2}+\left\vert J_{a+1}( 2Mr)
\right\vert ^{2}\\
=  &  \ \frac{2}{\pi Mr}\cos^{2}\left(  Mr-( a+1) \frac{\pi}{2}-\frac{\pi}%
{4}\right)  +\frac{2}{2\pi Mr}\cos^{2}\left(  2Mr-( a+1) \frac{\pi}{2}%
-\frac{\pi}{4}\right) \\
&  +O( m^{-2})\\
\geqslant &  \ \frac{C}{m}\left(  \cos^{2}\left(  Mr-( a+1) \frac{\pi}%
{2}-\frac{\pi}{4}\right)  +\cos^{2}\left(  2Mr-( a+1) \frac{\pi}{2}-\frac{\pi
}{4}\right)  \right)  +O( m^{-2}).
\end{align*}
Calling%
\[
\omega=Mr-( a+1) \frac{\pi}{2}-\frac{\pi}{4},
\]
we have that
\begin{align*}
&  \cos^{2}\left(  Mr-(a+1)\frac{\pi}{2}-\frac{\pi}{4}\right)  +\cos
^{2}\left(  2Mr-(a+1)\frac{\pi}{2}-\frac{\pi}{4}\right) \\
=  &  \ \cos^{2}(\omega)+\cos^{2}\left(  2\omega+(a+1)\frac{\pi}{2}+\frac{\pi
}{4}\right)  .
\end{align*}
We claim that this expression never vanishes. Assume the contrary, then
$\omega=\frac{\pi}{2}+k\pi$ for some integer $k,$ and $2\omega+(a+1)\frac{\pi
}{2}+\frac{\pi}{4}=\frac{\pi}{2}+h\pi$ for some integer $h$. This implies%
\[
\pi+2k\pi+(a+1)\frac{\pi}{2}+\frac{\pi}{4}=\frac{\pi}{2}+h\pi.
\]
Recalling that $a=\frac{d-2}{2}$, this identity implies that $\frac{d-1}{4}$
is an integer. This contradicts $d\not \equiv 1\,\,( \operatorname{mod}4) $.
Therefore for all $\omega\in\mathbb{R}$,%
\[
\cos^{2}( \omega) +\cos^{2}\left(  2\omega+(a+1)\frac{\pi}{2}+\frac{\pi}%
{4}\right)  \geqslant c,
\]
so that for large $m$%
\[
\left\vert J_{a+1}( Mr) \right\vert ^{2}+\left\vert J_{a+1}( 2Mr) \right\vert
^{2}\geqslant\frac{C}{m}.
\]
Finally,%
\begin{align*}
\left\vert \int_{B_{r}(o)}Z_{o}^{m}(x)d\mu(x)\right\vert ^{2}+\left\vert
\int_{B_{2r}(o)}Z_{o}^{m}(x)d\mu(x)\right\vert ^{2}  &  \geqslant
Cm^{-2a-3}d_{m}^{2}+d_{m}^{2}O( m^{-2a-4})\\
&  \geqslant Cd_{m}^{2}m^{-2a-3}.
\end{align*}

\end{proof}

We can now prove Theorem \ref{Thm 1}.

\begin{proof}
[Proof of Theorem \ref{Thm 1}]Let $m_{0}$ as in Lemma \ref{Lemma Due cerchi}
and let $X\geqslant m_{0}$. Then by Lemma \ref{Lemma Quadrato Discrep},
Proposition \ref{Prop Cassels} and Lemma \ref{Lemma Due cerchi}, we have%
\begin{align*}
&  \left\Vert D_{r}\right\Vert _{2}^{2}+\left\Vert D_{2r}\right\Vert _{2}%
^{2}\\
=  &  \ \sum_{m=1}^{+\infty} d_{m}^{-2}\left(  \left\vert \int_{B_{r}(o)}%
Z_{o}^{m}(x)d\mu(x)\right\vert ^{2}+\left\vert \int_{B_{2r}(o)}Z_{o}%
^{m}(x)d\mu(x)\right\vert ^{2}\right)  \sum_{\ell=1}^{d_{m}}\left\vert
\sum_{j=1}^{N}a_{j}Y_{m}^{\ell}(x_{j})\right\vert ^{2}\\
\geqslant &  \ \min_{m_{0}\leqslant m\leqslant X}\left(  d_{m}^{-2}\left(
\left\vert \int_{B_{r}(o)}Z_{o}^{m}(x)d\mu(x)\right\vert ^{2}+\left\vert
\int_{B_{2r}(o)}Z_{o}^{m}(x)d\mu(x)\right\vert ^{2}\right)  \right) \\
&  \times\sum_{m=m_{0}}^{X}\sum_{\ell=1}^{d_{m}}\left\vert \sum_{j=1}^{N}%
a_{j}Y_{m}^{\ell}(x_{j})\right\vert ^{2}\\
\geqslant &  \ C\left(  \min_{m_{0}\leqslant m\leqslant X}m^{-2a-3}\right)
\left(  C_{1}X^{d}\sum_{j=1}^{N}a_{j}^{2}-C_{0}\right)  \geqslant
CX^{-2a-3}\left(  C_{1}X^{d}\sum_{j=1}^{N}a_{j}^{2}-C_{0}\right)  .
\end{align*}
Applying Cauchy-Schwarz inequality to $\sum_{j=1}^{N}a_{j}=1$ gives
\[
\sum_{j=1}^{N}a_{j}^{2}\geqslant\frac{1}{N}.
\]
Let $N\geqslant N_{0}=m_{0}^{d}C_{1}/(2C_{0})$. Then, setting $X=\left[
(2C_{0}C_{1}^{-1}N)^{1/d}\right]  +1$, we have $X\geqslant m_{0}$, $C_{1}%
X^{d}\sum_{j=1}^{N}a_{j}^{2}-C_{0}\geqslant C_{0}$ and%
\begin{equation}
\left\Vert D_{r}\right\Vert _{2}^{2}+\left\Vert D_{2r}\right\Vert _{2}%
^{2}\geqslant CN^{-1-\frac{1}{d}}. \label{Discrep N}%
\end{equation}
Let now $N<N_{0}$ and let us consider the points and weights%
\[
\widetilde{x}_{j}=%
\begin{cases}
x_{j} & 1\leqslant j\leqslant N-1,\\
x_{N} & N\leqslant j\leqslant N_{0},
\end{cases}
\qquad\widetilde{a}_{j}=%
\begin{cases}
a_{j} & 1\leqslant j\leqslant N-1,\\
\dfrac{a_{N}}{N_{0}-N+1} & N\leqslant j\leqslant N_{0}.
\end{cases}
\]
Since the discrepancy $\widetilde{D}_{r}$ of the points $\{\widetilde{x}%
_{j}\}_{j=1}^{N_{0}}$ and weights $\{\widetilde{a}_{j}\}_{j=1}^{N_{0}}$
coincides with the discrepancy $D_{r}$ of the points $\{x_{j}\}_{j=1}^{N}$ and
weights $\{a_{j}\}_{j=1}^{N}$, applying (\ref{Discrep N}) to $\widetilde
{D}_{r}$ gives%
\[
\left\Vert D_{r}\right\Vert _{2}^{2}+\left\Vert D_{2r}\right\Vert _{2}%
^{2}\geqslant CN_{0}^{-1-\frac{1}{d}}=C\,\left(  \frac{N}{N_{0}}\right)
^{1+\frac{1}{d}}N^{-1-\frac{1}{d}}\geqslant C\left(  \frac{1}{N_{0}}\right)
^{1+\frac{1}{d}}N^{-1-\frac{1}{d}}%
\]
also when $1\leqslant N<N_{0}.$
\end{proof}

To prove Theorem \ref{Thm 2} we need the following result of Frenzen and Wong
(see \cite[Corollary 2]{F-W}) on the zeroes of the Jacobi polynomials .

\begin{theorem}
[Frenzen-Wong]\label{Thm Frenzen-Wong}Let $a\geqslant-\frac{3}{2}$,
$a+b\geqslant-3$, and let $0<\theta_{m-1,1}<\theta_{m-1,2}<\ldots
<\theta_{m-1,m-1}<\pi$ be the zeros of $P_{m-1}^{a+1,b+1}(\cos\theta)$. Then,
as $m\rightarrow+\infty$, we have%
\begin{align*}
\theta_{m-1,\ell}  &  =\frac{j_{a+1,\ell}}{M}+\frac{1}{M^{2}}\left\{  \left(
(a+1)^{2}-\frac{1}{4}\right)  \frac{1-t\cot t}{2t}-\frac{(a+1)^{2}-(b+1)^{2}%
}{4}\tan\frac{t}{2}\right\} \\
&  +t^{2}O\left(  \frac{1}{m^{3}}\right)
\end{align*}
where $j_{a+1,\ell}$ is the $\ell$-th positive zero of the Bessel function
$J_{a+1}(x)$ and $t=j_{a+1,\ell}/M$. The $O$-term is uniformly bounded for all
values of $\ell=1,\ldots,\left[  \gamma m\right]  $, where $\gamma\in(0,1)$.
\end{theorem}

\begin{lemma}
\label{Lemma MacMahon}Let $0<\gamma_{1}<\gamma_{2}<1$. For every $m\geqslant1$
and every $\left[  \gamma_{1}m\right]  \leqslant\ell\leqslant\left[
\gamma_{2}m\right]  $ we have%
\[
\theta_{m-1,\ell}=\frac{\ell\pi+a \frac{\pi}{2}+\frac{\pi}{4}}{M}+O\left(
\frac{1}{m^{2}}\right)  .
\]

\end{lemma}

\begin{proof}
This follows directly from McMahon's expansion (see e.g. (1.5) in \cite{Elb})%
\[
j_{a+1,\ell}=\ell\pi+ a \frac{\pi}{2}+\frac\pi4+O\left(  \frac{1}{\ell
}\right)
\]
and Theorem \ref{Thm Frenzen-Wong}.
\end{proof}

\begin{lemma}
\label{Jacobi da sotto}There exists $m_{0}>0$ such that for every $\delta>0$
and for almost every $r\in( 0,\pi) $ there exists a constant $C>0$ such that
for $m\geqslant m_{0}$,%
\[
\left\vert \sin^{2a+2}\left(  \frac{r}{2}\right)  \cos^{2b+2}\left(  \frac
{r}{2}\right)  P_{m-1}^{a+1,b+1}( \cos r) \right\vert \geqslant Cm^{-\frac
{3}{2}-\delta}.
\]

\end{lemma}

\begin{proof}
Let $\theta_{\ell}=\theta_{m-1,\ell}$, for $\ell=1,\ldots,m-1$, be the
positive zeros of $P_{m-1}^{a+1,b+1}(\cos\theta)$, and let%
\[
\mathcal{P}_{m-1}^{a+1,b+1}(r)=\sin^{2a+2}\left(  \frac{r}{2}\right)
\cos^{2b+2}\left(  \frac{r}{2}\right)  P_{m-1}^{a+1,b+1}(\cos r).
\]
Let $\varepsilon>0$. Let $\gamma_{1}$ and $\gamma_{2}$ be such that for $m$
sufficiently large%
\[
\frac{\varepsilon}{2}<\frac{[\gamma_{1}m]\pi+a\frac{\pi}{2}+\frac{\pi}{4}}%
{M}<\varepsilon
\]
and%
\[
\pi-\varepsilon<\frac{[\gamma_{2}m]\pi+a\frac{\pi}{2}+\frac{\pi}{4}}{M}%
<\pi-\frac{\varepsilon}{2}.
\]
Then, if $\eta\in(0,1)$,%
\[
\int_{\varepsilon}^{\pi-\varepsilon}\left\vert \mathcal{P}_{m-1}%
^{a+1,b+1}(r)\right\vert ^{-1+\eta}dr\leqslant\sum_{\ell=[\gamma_{1}%
m]}^{[\gamma_{2}m]-1}\int_{\theta_{\ell}}^{\theta_{\ell+1}}\left\vert
\mathcal{P}_{m-1}^{a+1,b+1}(r)\right\vert ^{-1+\eta}dr.
\]
By \eqref{rodrigues} we have%
\[
\frac{d}{dr}\mathcal{P}_{m-1}^{a+1,b+1}(r)=\frac12m\mathcal{P}_{m}%
^{a,b}(r)\sin(r).
\]
Hence, if $\varepsilon/2\leqslant r\leqslant\pi-\varepsilon/2$, then by
(\ref{Jacobi Bessel})%
\begin{align*}
\frac{d}{dr}\mathcal{P}_{m-1}^{a+1,b+1}(r)  &  =m\sin^{2a+1}\left(  \frac
{r}{2}\right)  \cos^{2b+1}\left(  \frac{r}{2}\right)  P_{m}^{a,b}(\cos r)\\
&  =m\sin^{a+1}\left(  \frac{r}{2}\right)  \cos^{b+1}\left(  \frac{r}%
{2}\right)  M^{-a}\frac{\Gamma(m+a+1)}{m!}\left(  \frac{r}{\sin r}\right)
^{\frac{1}{2}}J_{a}( Mr)\\
&  +\sin^{a+1}\left(  \frac{r}{2}\right)  \cos^{b+1}\left(  \frac{r}%
{2}\right)  r^{\frac{1}{2}}O(m^{-\frac{1}{2}}).
\end{align*}
Therefore%
\begin{equation}
\left\vert \frac{d}{dr}\mathcal{P}_{m-1}^{a+1,b+1}(r)\right\vert \geqslant
cm\left\vert J_{a}( Mr) \right\vert -O( m^{-\frac{1}{2}}) . \label{Derivata}%
\end{equation}
Let $\tau>0$ and let $r\in\left(  \theta_{\ell}-\frac{\tau}{M},\theta_{\ell
}+\frac{\tau}{M}\right)  $. Then by Lemma \ref{Lemma MacMahon}%
\[
Mr-a\frac{\pi}{2}-\frac{\pi}{4}\in\left(  \ell\pi+O\left(  \frac{1}{m}\right)
-\tau,\ell\pi+O\left(  \frac{1}{m}\right)  +\tau\right)  ,
\]
so that%
\[
\left\vert \cos\left(  Mr-a\frac{\pi}{2}-\frac{\pi}{4}\right)  \right\vert
\geqslant c.
\]
By the asymptotic expansion of the Bessel function \eqref{asymp Bessel} we
therefore obtain%
\begin{align*}
\left\vert J_{a}( Mr) \right\vert  &  \geqslant\sqrt{\frac{2}{\pi Mr}%
}\left\vert \cos\left(  Mr-a\frac{\pi}{2}-\frac{\pi}{4}\right)  \right\vert
-O( m^{-\frac{3}{2}}) \geqslant cm^{-\frac{1}{2}},
\end{align*}
and by (\ref{Derivata})%
\[
\left\vert \frac{d}{dr}\mathcal{P}_{m-1}^{a+1,b+1}(r)\right\vert \geqslant
cm^{\frac{1}{2}}.
\]
Thus,%
\begin{align*}
&  \int_{\theta_{\ell}}^{\theta_{\ell+1}}\left\vert \mathcal{P}_{m-1}%
^{a+1,b+1}(r)\right\vert ^{-1+\eta}dr\\
=  &  \int_{\theta_{\ell}}^{\theta_{\ell}+\frac{\tau}{M}}\left\vert
\mathcal{P}_{m-1}^{a+1,b+1}(r)\right\vert ^{-1+\eta}dr+\int_{\theta_{\ell
}+\frac{\tau}{M}}^{\theta_{\ell+1}-\frac{\tau}{M}}\left\vert \mathcal{P}%
_{m-1}^{a+1,b+1}(r)\right\vert ^{-1+\eta}dr\\
&  +\int_{\theta_{\ell+1}-\frac{\tau}{M}}^{\theta_{\ell+1}}\left\vert
\mathcal{P}_{m-1}^{a+1,b+1}(r)\right\vert ^{-1+\eta}dr\\
\leqslant &  \int_{\theta_{\ell}}^{\theta_{\ell}+\frac{\tau}{M}}\left\vert
cm^{1/2}(r-\theta_{\ell})\right\vert ^{-1+\eta}dr+\int_{\theta_{\ell}%
+\frac{\tau}{M}}^{\theta_{\ell+1}-\frac{\tau}{M}}\left\vert cm^{1/2}\left(
\frac{\tau}{M}\right)  \right\vert ^{-1+\eta}dr\\
&  +\int_{\theta_{\ell+1}-\frac{\tau}{M}}^{\theta_{\ell+1}}\left\vert
cm^{1/2}(\theta_{\ell+1}-r)\right\vert ^{-1+\eta}dr\\
\leqslant &  \ cm^{-\frac{1}{2}-\frac{\eta}{2}}.
\end{align*}
Therefore%
\begin{align*}
\int_{\varepsilon}^{\pi-\varepsilon}\left\vert \mathcal{P}_{m-1}%
^{a+1,b+1}(r)\right\vert ^{-1+\eta}dr  &  \leqslant\sum_{\ell=[\gamma_{1}%
m]}^{[\gamma_{2}m]-1}\int_{\theta_{\ell}}^{\theta_{\ell+1}}\left\vert
\mathcal{P}_{m-1}^{a+1,b+1}(r)\right\vert ^{-1+\eta}dr\\
&  \leqslant cm\cdot m^{-\frac{1}{2}-\frac{\eta}{2}}=cm^{\frac{1}{2}%
-\frac{\eta}{2}}.
\end{align*}
Finally%
\begin{align*}
\int_{\varepsilon}^{\pi-\varepsilon}\left(  \sum_{m=1}^{+\infty}m^{-\sigma
}\left\vert \mathcal{P}_{m-1}^{a+1,b+1}(r)\right\vert ^{-1+\eta}\right)  dr
&  =\sum_{m=1}^{+\infty}\int_{\varepsilon}^{\pi-\varepsilon}m^{-\sigma
}\left\vert \mathcal{P}_{m-1}^{a+1,b+1}(r)\right\vert ^{-1+\eta}dr\\
&  \leqslant c\sum_{m=1}^{+\infty}m^{\frac{1}{2}-\sigma-\frac{\eta}{2}}.
\end{align*}
If $\sigma>3/2-\eta/2$
\[
\sum_{m=1}^{+\infty}m^{-\sigma}\left\vert \mathcal{P}_{m-1}^{a+1,b+1}%
(r)\right\vert ^{-1+\eta}%
\]
converges for almost every $r\in(\varepsilon,\pi-\varepsilon)$ so that for
almost every $r\in(\varepsilon,\pi-\varepsilon)$ we have
\[
cm^{-\frac{\sigma}{1-\eta}}\leqslant\left\vert \mathcal{P}_{m-1}%
^{a+1,b+1}(r)\right\vert .
\]

\end{proof}

We are ready to prove Theorem \ref{Thm 2}.

\begin{proof}
[Proof of Theorem \ref{Thm 2}]Let $m_{0}$ be as in Lemma \ref{Jacobi da sotto}%
. By Lemma \ref{Lemma Quadrato Discrep}, Proposition \ref{Prop Cassels} and
Lemma \ref{Trasf bolla}, for every $X\geqslant m_{0}$ we have%
\begin{align*}
&  \int_{\mathcal{M}}\left\vert D_{r}(x)\right\vert ^{2}d\mu(x)\\
\geqslant &  \sum_{m=m_{0}}^{X}\sum_{\ell=1}^{d_{m}}\left\vert \sum_{j=1}%
^{N}a_{j}Y_{m}^{\ell}(x_{j})\right\vert ^{2}\min_{m_{0}\leqslant m\leqslant
X}d_{m}^{-2}\left\vert \int_{B_{r}(o)}Z_{o}^{m}(y)d\mu(y)\right\vert ^{2}\\
\geqslant &  \left(  C_{1}X^{d}\sum_{j=1}^{N}a_{j}^{2}-C_{0}\right) \\
&  \times\min_{m_{0}\leqslant m\leqslant X}\left\vert \frac{c(a,b)}%
{mP_{m}^{a,b}(1)}P_{m-1}^{a+1,b+1}(\cos r)\left(  \sin\left(  \frac{r}%
{2}\right)  \right)  ^{2a+2}\left(  \cos\left(  \frac{r}{2}\right)  \right)
^{2b+2}\right\vert ^{2}.
\end{align*}
Using Lemma \ref{Jacobi da sotto}, for almost every $r\in(0,\pi)$ and for
$m\geqslant m_{0}$ we obtain%
\begin{align*}
&  \min_{m_{0}\leqslant m\leqslant X}\left\vert \frac{c(a,b)}{mP_{m}^{a,b}%
(1)}P_{m-1}^{a+1,b+1}(\cos r)\left(  \sin\left(  \frac{r}{2}\right)  \right)
^{2a+2}\left(  \cos\left(  \frac{r}{2}\right)  \right)  ^{2b+2}\right\vert
^{2}\\
\geqslant &  \ c\min_{m_{0}\leqslant m\leqslant X}\left\vert \frac{1}%
{mP_{m}^{a,b}(1)}m^{-\frac{3}{2}-\delta}\right\vert ^{2}\geqslant c\min
_{m_{0}\leqslant m\leqslant X}\left\vert \frac{1}{m^{a+1}}m^{-\frac{3}%
{2}-\delta}\right\vert ^{2}\geqslant cX^{-3-2\delta-d}.
\end{align*}
Hence%
\begin{align*}
\int_{\mathcal{M}}|D_{r}(x)|^{2}d\mu(x)  &  \geqslant c\left(  C_{1}X^{d}%
\sum_{j=1}^{N}a_{j}^{2}-C_{0}\right)  X^{-3-2\delta-d}\\
&  \geqslant c\left(  C_{1}X^{d}\frac{1}{N}-C_{0}\right)  X^{-3-2\delta-d}.
\end{align*}
Let $N\geqslant N_{0}=m_{0}^{d}C_{1}/(2C_{0})$. Then, setting $X=\left[
(2C_{0}C_{1}^{-1}N)^{1/d}\right]  +1$, we have $X\geqslant m_{0}$, $C_{1}%
X^{d}N^{-1}-C_{0}\geqslant C_{0}$ and
\[
\int_{\mathcal{M}}|D_{r}(x)|^{2}d\mu(x)\geqslant cN^{-1-\frac{3}{d}%
-\frac{2\delta}{d}}%
\]
for all $N\geqslant N_{0}$. The same argument used in the proof of Theorem
\ref{Thm 1} gives the result for every $N\geqslant1$.
\end{proof}

\section{\label{Cubatures}Cubature formulas}

For a given $N\geqslant1$, let $\{a_{j}\}_{j=1}^{N}$ be a set of weights and
$\{x_{j}\}_{j=1}^{N}$ be a set of points in $\mathcal{M}$. We say that
$\left(  \{a_{j}\}_{j=1}^{N},\{x_{j}\}_{j=1}^{N}\right)  $ is a cubature of
strength $X$ on $\mathcal{M}$ if%
\[
\int_{\mathcal{M}}P(x)d\mu(x)=\sum_{j=1}^{N}a_{j}P(x_{j})
\]
for every%
\[
P\in%
{\displaystyle\bigoplus_{0\leqslant m\leqslant X}}
\mathcal{H}_{m}.
\]

In the following theorem we will show that under suitable assumptions a
cubature of strength $\kappa N^{1/d}$ has optimal discrepancy. This was
already observed by M. Skriganov in \cite[Corollary 2.1]{Skriganov}, although
in his result the discrepancy contains an integration in the radius $r$ of the
balls with respect to an absolutely continuous measure in $[ 0,\pi] $. The
same technique actually gives a stronger estimate, which is uniform in the
radius $r$.

\begin{theorem}
\label{Thm cubature alto}Let $A,B,\varepsilon$ and $\kappa$ be positive
constants. There exists $c>0$ such that if $N\geqslant1$, and $\left(
\{a_{j}\}_{j=1}^{N},\{x_{j}\}_{j=1}^{N}\right)  $ gives a cubature of strength
$\kappa N^{1/d}$ satisfying%
\begin{equation}
0\leqslant a_{j}\leqslant\dfrac{A}{N},\label{Condizione pesi}%
\end{equation}
and for every $x\in\mathcal{M}$%
\begin{equation}
\#\left\{  j:d(x_{j},x)\leqslant N^{-\frac{1}{d}}\right\}  \leqslant
B,\label{Ben separati}%
\end{equation}
then, for every $r\in\lbrack0,\pi-\varepsilon]$ we have
\begin{equation}
\int_{\mathcal{M}}|D_{r}(x)|^{2}d\mu(x)\leqslant cN^{-1-\frac{1}{d}%
}.\label{Stima dall'alto}%
\end{equation}

\end{theorem}

Observe that cubatures as required by the above theorem do exist. Indeed, in
\cite{EEGGT} it is proved that for every $A\geqslant1$ there exists a constant
$C(A,\mathcal{M})>0$, depending only on $A$ and $\mathcal{M}$, such that if
$N\geqslant C(A,\mathcal{M})X^{d}$ and the weights $\{a_{j}\}_{j=1}^{N}$
satisfy the conditions (\ref{Condizione pesi}) and%
\[
\sum_{j=1}^{N}a_{j}=1,
\]
then there is a choice of points $\{x_{j}\}_{j=1}^{N}$ such that $\left(
\{a_{j}\}_{j=1}^{N},\{x_{j}\}_{j=1}^{N}\right)  $ is a cubature of strength
$X$. The construction starts with a well separated set of points
$\{y_{j}\}_{j=1}^{N}$ and ends with a new set of points $\{x_{j}\}_{j=1}^{N}$
in such a way that each $x_{j}$ has distance at most $c_{1}(A,\mathcal{M}%
)N^{-1/d}$ from $y_{j}$. This implies that the number of points $x_{j}$
contained in any ball of radius $N^{-1/d}$ is uniformly bounded so that
(\ref{Ben separati}) is satisfied with a constant $B$ depending on
$\mathcal{M}$ and $A$. It suffices to apply the above construction with
$X=\kappa N^{1/d}$ ($\kappa\leqslant C(A,\mathcal{M})^{-1/d}$).

\begin{proof}
First of all we observe that it is enough to prove (\ref{Stima dall'alto}) for
$N$ sufficiently large. Let $N_{0}$ be a positive constant to be choosen later
and assume $N\geqslant N_{0}$. Since $\left(  \{a_{j}\}_{j=1}^{N}%
,\{x_{j}\}_{j=1}^{N}\right)  $ gives a cubature of strength $\kappa N^{1/d}$,
for every $m\leqslant\kappa N^{1/d}$ we have%
\[
\sum_{j=1}^{N}a_{j}Y_{m}^{\ell}(x_{j})=0.
\]
Then by Lemma \ref{Lemma Quadrato Discrep} and Lemma \ref{Trasf bolla} we have%
\begin{align*}
&  \int_{\mathcal{M}}|D(x)|^{2}d\mu(x)\\
= &  \sum_{m>\kappa N^{\frac{1}{d}}}\sum_{\ell=1}^{d_{m}}\left\vert \sum
_{j=1}^{N}a_{j}Y_{m}^{\ell}(x_{j})\right\vert ^{2}\left\vert c(a,b)\frac
{P_{m-1}^{a+1,b+1}(\cos r)}{mP_{m}^{a,b}(1)}\left(  \sin\left(  \frac{r}%
{2}\right)  \right)  ^{2a+2}\left(  \cos\left(  \frac{r}{2}\right)  \right)
^{2b+2}\right\vert ^{2}.
\end{align*}
We want to show that there exist $c>0$ and $\zeta>0$ such that for every
$r\in\lbrack0,\pi-\varepsilon]$, for every $N\geqslant1$ and for every
$m>\kappa N^{1/d}$%
\begin{align}
&  \left\vert P_{m-1}^{a+1,b+1}(\cos r)\left(  \sin\left(  \frac{r}{2}\right)
\right)  ^{2a+2}\left(  \cos\left(  \frac{r}{2}\right)  \right)
^{2b+2}\right\vert ^{2}\label{*}\\
\leqslant &  \ cN\int_{0}^{\zeta N^{-\frac{1}{d}}}\left\vert P_{m-1}%
^{a+1,b+1}(\cos u)\left(  \sin\left(  \frac{u}{2}\right)  \right)
^{2a+2}\left(  \cos\left(  \frac{u}{2}\right)  \right)  ^{2b+2}\right\vert
^{2}du.\nonumber
\end{align}
By (\ref{Jacobi Bessel}) we have the expansion%
\begin{align*}
&  \left(  \sin\frac{r}{2}\right)  ^{a+\frac{3}{2}}\left(  \cos\frac{r}%
{2}\right)  ^{b+\frac{3}{2}}P_{m-1}^{a+1,b+1}(\cos r)\\
= &  \ M^{-a-1}\frac{\Gamma(m+a+1)}{\sqrt{2}(m-1)!}\sqrt{r}J_{a+1}%
(Mr)+O(rm^{-\frac{3}{2}})
\end{align*}
uniformly in $r\in\lbrack0,\pi-\varepsilon]$ and $m>\kappa N^{\frac{1}{d}}$.
Hence%
\begin{align*}
&  \left\vert P_{m-1}^{a+1,b+1}(\cos r)\left(  \sin\left(  \frac{r}{2}\right)
\right)  ^{2a+2}\left(  \cos\left(  \frac{r}{2}\right)  \right)
^{2b+2}\right\vert ^{2}\\
\leqslant &  \left\vert P_{m-1}^{a+1,b+1}(\cos r)\left(  \sin\left(  \frac
{r}{2}\right)  \right)  ^{a+\frac{3}{2}}\left(  \cos\left(  \frac{r}%
{2}\right)  \right)  ^{b+\frac{3}{2}}\right\vert ^{2}\\
\leqslant &  \left\vert M^{-a-1}\frac{\Gamma(m+a+1)}{\sqrt{2}(m-1)!}\sqrt
{r}J_{a+1}(Mr)\right\vert ^{2}+O(rm^{-2})\\
\leqslant &  \ cm^{-1}+O(rm^{-2})\leqslant cm^{-1}.
\end{align*}
Let $0<\zeta<\pi N_{0}^{\frac{1}{d}}$. For every $m>\kappa N^{\frac{1}{d}}$,%
\begin{align*}
&  \int_{0}^{\zeta N^{-\frac{1}{d}}}\left\vert P_{m-1}^{a+1,b+1}(\cos
u)\left(  \sin\left(  \frac{u}{2}\right)  \right)  ^{2a+2}\left(  \cos\left(
\frac{u}{2}\right)  \right)  ^{2b+2}\right\vert ^{2}du\\
= &  \int_{0}^{\zeta N^{-\frac{1}{d}}}\left\vert P_{m-1}^{a+1,b+1}(\cos
u)\left(  \sin\left(  \frac{u}{2}\right)  \right)  ^{a+\frac{3}{2}}\left(
\cos\left(  \frac{u}{2}\right)  \right)  ^{b+\frac{3}{2}}\right\vert ^{2}\\
&  \times\left(  \sin\left(  \frac{u}{2}\right)  \right)  ^{2a+1}\left(
\cos\left(  \frac{u}{2}\right)  \right)  ^{2b+1}du\\
\geqslant &  c\int_{\zeta N^{-\frac{1}{d}}/2}^{\zeta N^{-\frac{1}{d}}%
}\left\vert P_{m-1}^{a+1,b+1}(\cos u)\left(  \sin\left(  \frac{u}{2}\right)
\right)  ^{a+\frac{3}{2}}\left(  \cos\left(  \frac{u}{2}\right)  \right)
^{b+\frac{3}{2}}\right\vert ^{2}u^{2a+1}du\\
\geqslant c &  \int_{\zeta N^{-\frac{1}{d}}/2}^{\zeta N^{-\frac{1}{d}}%
}\left\vert \sqrt{u}J_{a+1}(Mu)+O(um^{-\frac{3}{2}})\right\vert ^{2}%
u^{2a+1}du\\
\geqslant &  c\int_{\zeta N^{-\frac{1}{d}}/2}^{\zeta N^{-\frac{1}{d}}%
}\left\vert \sqrt{u}J_{a+1}(Mu)\right\vert ^{2}u^{2a+1}du-c\int_{\zeta
N^{-\frac{1}{d}}/2}^{\zeta N^{-\frac{1}{d}}}\left\vert \sqrt{u}J_{a+1}%
(Mu)O(um^{-\frac{3}{2}})\right\vert u^{2a+1}du\\
= &  \mathcal{A}-\mathcal{B}.
\end{align*}
To estimate $\mathcal{A}$ we use (\ref{asymp Bessel}). We have
\begin{align*}
\mathcal{A} &  =c\int_{\zeta N^{-\frac{1}{d}}/2}^{\zeta N^{-\frac{1}{d}}%
}\left\vert \sqrt{u}J_{a+1}(Mu)\right\vert ^{2}u^{2a+1}du\\
&  =c\int_{\zeta N^{-\frac{1}{d}}/2}^{\zeta N^{-\frac{1}{d}}}\left\vert
\left(  M^{-\frac{1}{2}}\cos\left(  Mu-(a+1)\frac{\pi}{2}-\frac{\pi}%
{4}\right)  +O(m^{-\frac{3}{2}}u^{-1})\right)  \right\vert ^{2}u^{2a+1}du\\
&  \geqslant cm^{-1}\int_{\zeta N^{-\frac{1}{d}}/2}^{\zeta N^{-\frac{1}{d}}%
}\left\vert \cos\left(  Mu-(a+1)\frac{\pi}{2}-\frac{\pi}{4}\right)
\right\vert ^{2}u^{d-1}du-c\int_{\zeta N^{-\frac{1}{d}}/2}^{\zeta N^{-\frac
{1}{d}}}m^{-2}u^{-1}u^{2a+1}du\\
&  \geqslant cm^{-1}(\zeta N^{-\frac{1}{d}})^{d-1}\int_{\zeta N^{-\frac{1}{d}%
}/2}^{\zeta N^{-\frac{1}{d}}}\left\vert \cos\left(  Mu-(a+1)\frac{\pi}%
{2}-\frac{\pi}{4}\right)  \right\vert ^{2}du-cm^{-2}(\zeta N^{-\frac{1}{d}%
})^{d-1}\\
&  \geqslant cm^{-2}(\zeta N^{-\frac{1}{d}})^{d-1}\int_{M\zeta N^{-\frac{1}%
{d}}/2}^{M\zeta N^{-\frac{1}{d}}}\left\vert \cos\left(  v-(a+1)\frac{\pi}%
{2}-\frac{\pi}{4}\right)  \right\vert ^{2}dv-cm^{-2}(\zeta N^{-\frac{1}{d}%
})^{d-1}\\
&  \geqslant cm^{-2}(\zeta N^{-\frac{1}{d}})^{d-1}M\zeta N^{-\frac{1}{d}%
}-cm^{-2}(\zeta N^{-\frac{1}{d}})^{d-1}\\
&  \geqslant cm^{-1}\zeta^{d-1}N^{-1}(\zeta-cm^{-1}N^{\frac{1}{d}})\geqslant
cm^{-1}\zeta^{d-1}N^{-1}(\zeta-c\kappa^{-1})\geqslant cm^{-1}\zeta^{d-1}N^{-1}%
\end{align*}
for $\zeta>2\kappa^{-1}c$ (and therefore for $N_{0}$ large enough). Now we
estimate $\mathcal{B}$:%
\begin{align*}
\mathcal{B} &  =c\int_{\zeta N^{-\frac{1}{d}}/2}^{\zeta N^{-\frac{1}{d}}%
}\left\vert \sqrt{u}J_{a+1}(Mu)O(um^{-\frac{3}{2}})\right\vert u^{2a+1}du\\
&  \leqslant c\int_{\zeta N^{-\frac{1}{d}}/2}^{\zeta N^{-\frac{1}{d}}}%
m^{-2}u^{d}du\leqslant cm^{-2}(\zeta N^{-\frac{1}{d}})^{d+1}=cm^{-2}%
\zeta^{d+1}N^{-1-\frac{1}{d}}.
\end{align*}
It follows that%
\begin{align*}
&  \int_{0}^{\zeta N^{-\frac{1}{d}}}\left\vert P_{m-1}^{a+1,b+1}(\cos
u)\left(  \sin\left(  \frac{u}{2}\right)  \right)  ^{2a+2}\left(  \cos\left(
\frac{u}{2}\right)  \right)  ^{2b+2}\right\vert ^{2}du\\
\geqslant &  \ cm^{-1}\zeta^{d-1}N^{-1}-cm^{-2}\zeta^{d+1}N^{-1-\frac{1}{d}%
}=c\zeta^{d-1}m^{-1}N^{-1}\left(  1-cm^{-1}\zeta^{2}N^{-\frac{1}{d}}\right)
\\
\geqslant &  \ c\zeta^{d-1}m^{-1}N^{-1}\left(  1-cN^{-\frac{2}{d}}\right)
\geqslant cm^{-1}N^{-1},
\end{align*}
and (\ref{*}) is proved. Thus,%
\begin{align*}
&  \int_{\mathcal{M}}|D_{r}(x)|^{2}d\mu(x)\\
= &  \sum_{m>\kappa N^{\frac{1}{d}}}\sum_{\ell=1}^{d_{m}}\left\vert \sum
_{j=1}^{N}a_{j}Y_{m}^{\ell}(x_{j})\right\vert ^{2}\left\vert c(a,b)\frac
{P_{m-1}^{a+1,b+1}(\cos r)}{mP_{m}^{a,b}(1)}\left(  \sin\left(  \frac{r}%
{2}\right)  \right)  ^{2a+2}\left(  \cos\left(  \frac{r}{2}\right)  \right)
^{2b+2}\right\vert ^{2}\\
\leqslant &  \ cN\sum_{m>\kappa N^{\frac{1}{d}}}\sum_{\ell=1}^{d_{m}%
}\left\vert \sum_{j=1}^{N}a_{j}Y_{m}^{\ell}(x_{j})\right\vert ^{2}\\
&  \times\int_{0}^{\zeta N^{-\frac{1}{d}}}\left\vert c(a,b)\frac
{P_{m-1}^{a+1,b+1}(\cos u)}{mP_{m}^{a,b}(1)}\left(  \sin\left(  \frac{u}%
{2}\right)  \right)  ^{2a+2}\left(  \cos\left(  \frac{u}{2}\right)  \right)
^{2b+2}\right\vert ^{2}du\\
= &  \ cN\int_{0}^{\zeta N^{-\frac{1}{d}}}\int_{\mathcal{M}}\left\vert
D_{u}(x)\right\vert ^{2}d\mu(x)du.
\end{align*}
Since $u\leqslant\zeta N^{-1/d}$ using (\ref{Ben separati}) and noticing that
by Lemma \ref{Lemma misura} a ball of radius $u$ can be covered by
$\frac{c_{2}(d,d_{0})}{c_{1}(d,d_{0})}(2(\zeta+1))^{d}$ balls of radius
$N^{-1/d}$ we have%
\begin{align*}
\left\vert D_{u}(x)\right\vert  &  =\left\vert \sum_{j=1}^{N}a_{j}\chi
_{B_{u}(x)}(x_{j})-\mu(B_{u}(x))\right\vert \\
&  \leqslant\frac{A}{N}\sum_{j=1}^{N}\chi_{B_{u}(x)}(x_{j})+\mu(B_{u}(x))\\
&  \leqslant\frac{A}{N}\frac{c_{2}(d,d_{0})}{c_{1}(d,d_{0})}(2(\zeta
+1))^{d}B+\frac{c_{2}(d,d_{0})\zeta^{d}}{N}\leqslant\frac{c}{N},
\end{align*}
so that%
\[
\int_{\mathcal{M}}|D_{u}(x)|^{2}d\mu(x)\leqslant\frac{C}{N^{2}}.
\]
Finally,%
\begin{align*}
\int_{\mathcal{M}}\left\vert D_{r}(x)\right\vert ^{2}d\mu(x) &  \leqslant
cN\int_{0}^{\zeta N^{-\frac{1}{d}}}\int_{\mathcal{M}}\left\vert D_{u}%
(x)\right\vert ^{2}d\mu(x)du\\
&  \leqslant cN\int_{0}^{\zeta N^{-\frac{1}{d}}}\frac{C}{N^{2}}du=cN^{-1-\frac
{1}{d}}.
\end{align*}

\end{proof}

\end{document}